\documentclass[11pt]{amsart}

\usepackage{amssymb,amsfonts,amscd}
\newtheorem{thm}{Theorem}[section]

\newtheorem{prop}[thm]{Proposition}

\newtheorem{mthm}[thm]{Main Theorem}

\theoremstyle{definition}
\newtheorem*{ex}{Example}

\newtheorem{defi}[thm]{Definition}

\newtheorem{rem}{Remark}

\author[E.,W.,Graczy\'nska]{Ewa Wanda Graczy\'nska}
\address{Department of Mathematics\\ Technical University, 45-036 Opole,
\\ ul. Luboszycka 3, POLAND}
\email{egracz@po.opole.pl}

\title[On the problem of basis for hyperquasivarieties]{On the problem of basis for hyperquasivarieties}

\thanks{This work was partially supported by Technische Universit\"{a}t Dresden, Germany.}

\usepackage{aaa68}
\AAAsubmitted{July 20, 2004}
\AAApage{1}
\begin{document}
\begin{abstract}
Our aim is to present a solution of the Problem 31 of Chapter 8 of \cite{2}
by presenting a general theorem on hyper-quasi-identities considered as
hyperbasis for solid quasi-varieties. In the last section we point out also
a generalization of that result for the notion of M-hyper-quasi-identities.

The results were presented at AAA68 Conference at
Technische Universit\"{a}t Dresden, June 13, 2004.
\end{abstract}
\maketitle

\section{Notations}

An identity is a pair of terms where the variables are bound by
universal quantifiers. Let us take the following medial identity as
an example
\[
\forall u\forall x\forall y\forall w\,(u\cdot x)\cdot (y\cdot
w)=(u\cdot y)\cdot (x\cdot w).
\]
Let us look at the following hyperidentity
\[
\,\forall f\forall u\forall x\forall y\forall
w\,f(f(u,v),f(x,y))=f(f(u,x),f(v,y)).
\]
The hypervariable $f$ is considered in a very specific way. Firstly
every hypervariable is restricted to functional symbols of a given
arity. Secondly, for a hypervariable only functions of a given arity
have to be substituted, and these functions have to be term
functions. Let us take the variety $A_{n,0}$ of abelian groups of
finite exponent $n.$ Every binary term $t\equiv t(x,y)$ can be
presented by $t(x,y)=ax+by$ with $a,b\in \mathbb N_{0}.$ If we
substitute the binary hypervariable $F$ in the above hyperidentity
by $ax+by,$ leaving its variables unchanged, we get
\[
a(au+bv)+b(ax+by)=a(au+bx)+b(av+by).
\]
This identity holds for every term $t(x,y)=ax+by$ for the variety
$A_{n,0}$. Therefore we say that the hyperidentity holds for the
variety $A_{n,0}.$

\section {Hyper-quasi-identities}

In the sequel we avoid to use the definition of {\em hyperterm} from
\cite{4} and accept the convention that a {\em hyperterm} is
formally the same as a {\em term}. We accept the notation of
\cite{PMC}, \cite{5}, \cite{VAG}.

{\em Hypersubstitutions} of terms were defined in \cite{1}, \cite{2}
and \cite{4}. Shortly speaking, they are mappings sending terms to
terms by substituting variables by (the same) variables and
functional symbols by terms of the same arities, i.e. $\sigma(x) =
x$ for any variable $x$, and for given functional symbol $f$, assume
that $\sigma(f)$ is a given term of the same arity as $f$, then
$\sigma$ acts on all terms of a given type in an inductive way:
$\sigma(f(p_1,...,p_n)) = \sigma(f)(\sigma(p_1),...,\sigma(p_n))$.

$\Sigma(\tau)$ denotes the set of all hypersubstitutions $\sigma$ of
a given type $\tau$.

We recall only our definitions of \cite{4} of hyperidentities
satisfied in an algebra of a given type and the notion of a {\em
hypervariety}:

\begin{defi}\label{D:2.1}
An algebra ${\bf A}$ {\it satisfies a hyperidentity} $h_{1}=h_{2}$ if for
every substitution of the hypervariables by terms (of the same
arity) of ${\bf A}$ leaving the variables unchanged, the identities
which arise hold in ${\bf A}$. In this case, we write ${\bf A}
\models (h_{1}=h_{2})$. A variety $V$ satisfies a hyperidentity
$h_{1}=h_{2}$ if every algebra in the variety does.
\end{defi}
\begin{defi}\label{D:2.2}
A class $V$ of a algebras of a given type is called a {\it hypervariety}
if and only if $V$ is defined by a set of hyperidentities.
\end{defi}
\begin{rem}
Some authors avoid to use our concepts since in \cite{WT} W. Taylor
defined hypervarieties to be classes of varieties of different types
satisfying certain sets of equations as identities (see also
\cite{WN}).
\end{rem}
The following was proved in \cite{4}:

\begin{thm} \label{T:2.3}
A variety $V$ of type $\tau$ is defined by a set of hyperidentities
if and only if $V = HSPD(V)$, i.e. $V$ is a variety closed under
derived algebras of type $\tau$.
\end{thm}
We recall from \cite{AIM} and \cite{VAG}:
\begin{defi}\label{D:2.4}
A {\it quasi-identity} $e$ is an implication of the form:
\begin{center}
(1.1) $(t_{0} = s_{0}) \wedge ... \wedge (t_{n-1} = s_{n-1})
\rightarrow (t_{n} = s_{n})$.
\end{center}
where $t_{i} = s_{i}$ are $k$-ary identities of a given type, for $i
= 0,...,n$.

A {\em quasi-identity} above is {\em satisfied in an algebra} ${\bf
A}$ of a given type if and only if the following implication is
satisfied in ${\bf A}$:  given a sequence $a_{1},...,a_{k}$ of
elements of $A$. If these elements satisfy the equations
$t_{i}(a_{1},...,a_{k}) = s_{i}(a_{1},...,a_{n})$ in ${\bf A}$, for
$i = 0, 1,..., n-1$, then the equality $t_{n}(a_{1},...,a_{k}) =
s_{n}(a_{1},...,a_{k})$ is satisfied in ${\bf A}$. In that case we
write:
\begin{center}
${\bf A} \models (t_{0} = s_{0}) \wedge ... \wedge (t_{n-1} =
s_{n-1}) \rightarrow (t_{n} = s_{n})$.
\end{center}

A quasi-identity $e$ is {\em satisfied in a class} $V$ of algebras
of a given type, if and only if it is satisfied in all algebras
${\bf A}$ belonging to $V$.

\end{defi}

Following A. I. Mal'cev \cite{AIM} we consider classes $QV$ of
algebras A of a given type $\tau$ defined by quasi-identities and
call them {\em quasivarieties}.

Formally (cf. \cite{2}, \cite{4}), a {\em hyper-quasi-identity} $e$
is the same as quasi-identity. Following the ideas of \cite{9}, part
5, cf. \cite{2} and \cite{EGDS2} (cf. \cite{CCKD}) we modify the
definition above in the following way.

\begin{defi}
A hyper-quasi-identity $e$ is {\it satisfied} (is {\it hyper-satisfied},
{\it holds})
in an algebra ${\bf A}$ if and only if the following implication is satisfied:\\
if $\sigma$ is a hypersubstitution of type $\tau$ and the elements
$a_{1},...,a_{n} \in A$ satisfy the equalities
$\sigma(t_{i})(a_{1},...,a_{k}) = \sigma(s_{i})(a_{1},...,a_{k})$ in
${\bf A}$, for $n = 0,1,...,n-1$, then the equality
$\sigma(t_{n})(a_{1},...,a_{k}) = \sigma(s_{n})(a_{1},...,a_{k})$
holds in ${\bf A}$.
\begin{center}
In that case, we write $V \models^H e$.

\end{center}
\end{defi}

In other words, hyper-quasi-identity is a universally closed  Horn
$\forall x \forall \sigma$-formulas, where x varies over all
sequences of individual variables (occurring in terms of the
implication) and $\sigma$ varies over all hypersubstitutions of a
given type. Our modification coincides with Definition 5.1.3 of
\cite{9} (cf. Definition 2.3 of \cite{CCKD}).

\begin{rem}
All hyper-quasi-identities and hyperidentities are written without
quantifiers but they are considered as universally closed Horn
$\forall$-formulas (cf. \cite{AIM}). A syntactic side of the notions
described here is  considered in a forthcoming paper, together with
a suitable generalization of the notion of hyperquasivariety to
M-hyperquasivariety.
\end{rem}

Let $V$ be a class of algebras of type $\tau$. Derived algebras were
defined in \cite{PMC}. Derived algebras of a given type $\tau$ were
defined in \cite{9}, \cite{4}.
\begin{defi}
Let  ${\bf A} =(A,F)$ be an algebra in $V$ and $\sigma$ a hypersubstitution
in $\Sigma(\tau)$. Then the algebra ${\bf B} = (A, (F)^\sigma)$ is a
{\it derived algebra} of ${\bf A}$, with the same universe $A$ and the set
$(F)^\sigma$ of all derived operations of $F$. We denote then ${\bf
B}$ as ${\bf A^\sigma}$.
\end{defi}

\section{Problem Formulation}

In the Example on p. 155 of \cite{2} the
authors proved  that the quasi-identities of the given basis of
\cite{2} are hyper-quasi-identities. In the Problem 31 of Chapter 8
they posed the question:

(3.1) {\em Prove that it is sufficient to prove that the quasi-identities of the
given basis of Example on p. 155 of \cite{2} are hyper-quasi-identities
for checking that every quasi-identity is a hyperquasi-identity}.

We shall refer to the question above as to the Problem (3.1).

In \cite{2} {\em solid quasivarieties} were defined as
quasivarieties closed under taking of derived algebras, i.e. we
assume the following
\begin{defi}
Let $QV$ be a quasivariety, then $QV$ is {\it solid} if and only if every
derived algebra ${\bf A^\sigma}$ belongs to $QV$, for every algebra
${\bf A}$ in $QV$ and $\sigma$ in $\Sigma(\tau)$.
\end{defi}
We write then, that
\begin{center}
$QV = D(QV)$
\end{center}
In \cite{EGDS2} {\em solid quasivarieties} in the above sense were
called {\em hyperquasivarieties}.

\subsection{Basis and hyperbasis}

Let $\Sigma$ be a set of quasi-identities of type $\tau$.
\begin{defi}
A quasivariety $QV$ of all algebras of type $\tau$ satisfying all
quasi-identities of $\Sigma$ is called a {\it quasivariety defined by
the basis} $\Sigma$.
\end{defi}

\begin{defi}
A quasivariety $QV$ of all algebras of type $\tau$ satisfying all
quasi-identities of $\Sigma$ as hyperquasi-identities is called a
{\it quasivariety defined by the hyperbasis} $\Sigma$.
\end{defi}

Recall Example of \cite{2}, p. 155:
\begin{ex}
Consider the quasivariety $QV$ of type $\tau=(2)$ defined by the following
identities:

(S1) $x(yz) = (xy)z$, (identities are regarded as quasi-identities),

(S2) $xx = x$,

(S3) $(xy)(uv) = (xu)(yv)$,

(S4) $(xy = yx) \rightarrow (x=y)$.
\end{ex}

The example above shows, that the given basis (S1)-(S4) is also a hyperbasis
of $QV$.

The authors of \cite{2} proved that the identities (S1), (S2), (S3)
are satisfied as hyperidentities and thus as hyperquasi-identities
in $QV$. Moreover (S4) is satisfied as a hyperquasi-identity.
Finally they concluded that the quasivariety $QV$ is solid.

Let us note, that in the proof in \cite{2} on p. 155 of the fact above,
the authors did not use their own definition on p. 155 of a {\it solid
qusivariety}. Therefore we present here an explicit proof of
the positive solution of the Problem 3.1 to make the situation more clear
for the reader. Namely, our main theorem states that every quasivariety
defined by a hyperbasis is solid and vice versa.

\begin{rem}
The positive answer of the Problem 3.1  can also be concluded via
Theorem 2.6 of \cite{CCKD}, even for the case of M-hyperquasi-equational
theories. For varieties the solution can be also
proved via Theorems 13.3 - 13.5 of \cite{EG1} or Theorem 14.34 of \cite{3}.
In a more general setting, the solution is also implicitly contained in
\cite{3}, Theorem 13.1.6 and Theorem 13.1.5 as a part of the  theory
of conjugate pairs of additive closure operators (cf. Lemma 4.9.2 and Theorem
4.9.3 of \cite{2}, p. 156). However, in the paper we do not use this theory,
neither the authors of \cite{CCKD} and  \cite{3} wrote
an explicit proof of the problem.

In addition, our Propositions 5.2 and 5.4 show that the Problem 3.1 considered
here is mainly connected with  relations between rules of inferences and not
directly with the quoted above results of \cite{CCKD} and \cite{2}, \cite{3}.
\end{rem}

\section{Problem Solution}
We present here a short proof of of the required statement
by proving the following:

\begin{mthm}
Let $QV$ be a  quasivariety  defined by a set $\Sigma$ of
quasi-identities. These quasi-identities are satisfied in $QV$ as
hyper-quasi-identities if and only if $QV$ is solid (is a
hyperquasivariety). Moreover, each quasi-identity satisfied in a
solid quasivariety $QV$ is satisfied in $QV$ as a
hyper-quasi-identity and vice versa.
\end{mthm}
\begin{proof}
Let be given  a quasivariety $QV$ with a basis $\Sigma$ of
quasi-identities. Assume that all quasi-identities of $\Sigma$  are
satisfied in $QV$ as hyper-quasi-identities (i.e. $\Sigma$ is a
hyperbasis of $QV$). Let ${\bf A}$ be an algebra in $QV$ and a
derived algebra ${\bf B} = A^\sigma$. Then for any quasi-identity
$e$ of  $\Sigma$, this quasi-identity $e$ is satisfied in ${\bf B}$,
for the following reason: each term $t$ of type $\tau$ is realized
in ${\bf B}$ as $\sigma(t)$, therefore for a given sequence
$a_1,...,a_k$ of elements of $A$, if these elements satisfy the
equations $t_i(a_1,...,a_k) = s_i(a_1,...,a_k)$ for $i = 1,...,n-1$,
in ${\bf B}$ (i.e. equations $\sigma(t_i)(a_1,...,a_k) =
\sigma(s_i)(a_1,...,a_k)$ for $i = 1,...,n-1$ in $A$) then the
equality $t_n(a_1,...,a_k)= s_n(a_1,...,a_k)$ is satisfied in ${\bf
B}$ (i.e. the equality $\sigma(t_n)(a_1,...,a_k)=
\sigma(s_n)(a_1,...,a_k)$ is satisfied in ${\bf A}$) from the
assumption that $\Sigma$ is a hyperbasis of $QV$. Therefore ${\bf
B}$ is in $QV$ and $QV$ is solid.

Let $QV$ be a  solid quasivariety  defined by a basis
$\Sigma$. Then any quasi-identity $e$ of $\Sigma$ is satisfied in
$QV$ as a hyper-quasi-identity (i.e.  $\Sigma$ is a hyperbasis of
$QV$). To show this take a $\sigma$ of  $\Sigma(\tau)$ and any
algebra ${\bf A}$ of $QV$. Assume that for a given sequence
$a_1,...,a_k $ of elements of $A$, these elements satisfy the
equations $\sigma(t_i)(a_1,...,a_k) = \sigma(s_i)(a_1,...,a_k)$ for
$i = 1,...,n-1$, then the equality $\sigma(t_n)(a_1,...,a_k)=
\sigma(s_n)(a_1,...,a_k)$ is satisfied in the derived algebra ${\bf
B} = {\bf A^\sigma}$, by similar arguments as above, as ${\bf B}$ is
in $QV$ by the assumption that $QV$ is solid and $e$ is satisfied in
$QV$ as an element of $\Sigma$.

Therefore we have proved that a quasivariety $QV$ defined by a basis
$\Sigma$ is solid if and only if $\Sigma$ is its hyperbasis. Let us
note, that every quasivariety $QV$ is defined by a basis $\Sigma$
consisting of all quasi-identities satisfied in $QV$. This
observation proves the last statement of  our main theorem.
\end{proof}

\section{Varieties}
As a specific case of {\em quasivarieties} one may consider
{\em varieties} of algebras and conclude a similar theorem for
{\em basis} and {\em hyperbasis of identities} (cf. \cite{4}):

\begin{prop}
Let $V$ be  a variety defined by a set $\Sigma$ of identities. These
identities are satisfied in $V$ as hyper-identities if and only if
$V$ is solid. Moreover, each identity satisfied in a solid variety
$V$ is satisfied in $V$ as a hyper-identity, i.e. $V$ is a
hypervariety and vice versa. Therefore the both notion coincide.
\end{prop}
\begin{proof}
The proof follows immediately from the Birkhoff's type theorem
proved in \cite{4}: a variety $V$ is solid if and only if
\begin{center}
$V = HSPD(V)$.
\end{center}
\end{proof}

Let us denote by $E$, the closure operator defined by the classical
rules (1)-(5) of inferences for identities (cf. \cite{GB}, \cite{5})
and $E^{H}$ denotes the closure operator defined by the rules (1) - (5)
and so called {\it hypersubstitution rule} (6) of \cite{4}, p. 308.
In \cite{EG} we proved the following Lemma 2, p. 121 (cf. Remark 1.2
of \cite{4} p. 308):
\begin{prop}
If $\Sigma$ is a set of identities of type $\tau$, closed under
hypersubstitution rule (6), then $E(\Sigma) = E^{H}(\Sigma)$.
\end{prop}  Let us note that the Proposition 5.2 can be slightly generalized by
considering the rule $(6)_{M}$, i.e. the M-hypersubstitution rule of
\cite{EG1}, p. 90, instead of the rule (6) (i.e. by considering only
hypersubstitutions from a given monoid $M$ instead all
hypersubstitutions of a given type). Let $E^{H}_{M}$ denotes the
closure operator defined by the inference rules (1) - (5) of G.
Birkhoff and the rule $(6)_{M}$. First recall from \cite{EG1}, p.
90:
\begin{defi}
Let be given a monoid $M$ of hypersubstitions of type $\tau$. The
{\it M-hypersubstition rule} is defined by the following:
\begin{center}
$(6)_{M}$ from $p = q$ conclude $\sigma(p) = \sigma(q)$, for every $\sigma \in M$.
\end{center}

\end{defi}
We get:
\begin{prop}
If $\Sigma$ is a set of identities of type $\tau$, closed under
M-hyper\-substitution rule, $(6)_{M}$, then $E(\Sigma) = E^{H}_{M}(\Sigma)$.
\end{prop}
\begin{proof} A proof is similar as those of \cite{EG}, p. 121-122.
The only difference is that considered hypersubstituions $\sigma$
may vary only within members of $M$.
\end{proof}

Recall from \cite{2}  that hyperidentities (especially in the theory
of Boolean algebras) has natural interpretations in the theory of
{\em swiching circuts}. Therefore our observation may have some
applications in algebraic computation.

\section{Conclusion} By other words, we observed that solid
quasivarieties are always defined by a hyperbasis. Moreover, the set
of all quasi-identities satisfied in a quasivariety $QV$ is
satisfied in $QV$ as hyperquasi-dentities if and only if $QV$ is
solid.

A suitable generalization of our observations made for the set
$\Sigma(\tau)$ of all hypersubstitutions can be extended to any
subset of $\Sigma(\tau)$, closed under superposition. This
generalization gives rise to so called {\em M-hypersubstitutions} of
a given type. We consider the set $\Sigma(\tau)$ of all
hypersubstitutions of a given type with the superposition operation
$\circ$ i.e. with  the composition of hypersubstitutions) and
consider the structure $(\Sigma(\tau), \circ)$ as a {\em monoid}. We
call it the {\em monoid of all hypersubstitutions of a given type
$\tau$}.

Let $M$ be a subset of $\Sigma(\tau)$ closed under superposition
$\circ$, i.e. a submonoid $(M, \circ)$ of the monoid $(\Sigma(\tau),
\circ)$ .
\begin{defi}
A (hyper)quasi-identity $e$ is satisfied as an {\it
M-hyper-quasi-identity} ({\it is M-hypersatisfied}) in an algebra
${\bf A}$ if and only if the following implication is satisfied:

if $\sigma$ is a hypersubstitution of type $\tau$ from the set $M$
and the elements $a_1,...,a_k$ of $A$ satisfy the equalities:
$\sigma(t_i)(a_1,...,a_k) =  \sigma(s_i)(a_1,...,a_k)$  in $A$ for
$i = 1,...,n-1$, then the equality $\sigma(t_n)(a_1,...,a_k) =
\sigma (s_n)(a_1,...,a_k)$ holds in $A$.
\end{defi}
In that case we write:
\begin{center}
${\bf A} \models_M^H (t_{0} = s_{0}) \wedge ... \wedge (t_{n-1} =
s_{n-1}) \rightarrow (t_{n} = s_{n})$.
\end{center}

\begin{defi}
Let $QV$ be a {\it a quasivariety}, then $QV$ is {\it M-solid} if and only if
every M-derived algebra ${\bf A^\sigma}$ belongs to $QV$, for every
algebra ${\bf A}$ in $QV$ and $\sigma$ in $M$.
\end{defi}
We write then, that
\begin{center}
$QV = D_M(QV)$
\end{center}
In \cite{EGDS3} M-solid quasivarieties were called M-hyperquasivarieties.

We finalize our considerations by formulating a generalization of
our main theorem for the case of {\em M-hypersubstitutions} without
giving details:
\begin{mthm}
Let $QV$ be a quasivariety  defined by a set $\Sigma$ of
quasi-identities. These quasi-identities are satisfied in $QV$ as an
M-hyper-quasi-identities if and only if $QV$ is M-solid (is an
M-hyperquasivariety). Moreover, each quasi-identity satisfied in an
M-solid quasivariety $QV$ is satisfied in $QV$ as an
M-hyper-quasi-identity and vice versa.
\end{mthm}
\begin{proof}
A proof is similar as those of the Main Theorem 4.1. The only
restriction is that hypersubstitions $\sigma$ and derived algebras
${\bf B} = {\bf A^\sigma}$ are considered for all $\sigma$ from a
given monoid $M$.
\end{proof}
\begin{prop}
Let $V$ be a variety defined by a set $\Sigma$ of identities. These
identities are satisfied in $V$ as M-hyper-identities if and only if
$V$ is M-solid. Moreover, each identity satisfied in an M-solid
variety $V$ is satisfied in $V$ as an M-hyper-identity, i.e. $V$ is
an M-hypervariety and vice versa.
\end{prop}
\begin{proof}
The proof follows immediately from the Birkhoff's type theorem
proved in \cite{EG1}: a variety $V$ is M-solid if and only if
\begin{center}
$V = HSPD_M(V)$.
\end{center}
\end{proof}
\begin{rem}
Let us note that in case $M$ is a trivial (i.e. 1-element) monoid of
hypersubstitutions of a given type $\tau$, then the satisfaction
$\models_M^H$ gives rise to the satisfaction $\models$ and the
operator $D_M$ to the  identity operator.

In case $M = \Sigma(\tau)$ we get the notion of $\models^H$
considered in \cite{EGDS2}.
\end{rem}

{\bf Acknowledgements}

The author expresses her thanks to the referees for their valuable comments.

\end{document}